\documentclass[11pt]{article}
\usepackage[utf8]{inputenc}
\usepackage[a4paper,margin=1in]{geometry}
\usepackage{amsmath,amssymb}
\usepackage{graphicx}
\usepackage{booktabs}
\usepackage{caption}
\usepackage{authblk}
\usepackage[hidelinks]{hyperref}
\title{Flocking Time and Graph Volume in Nonlinear Cucker--Smale Models with Dynamic Networks}
\author[1]{Sooyeon Song}
\author[1]{Jea-Hyun Park\thanks{Corresponding author. E-mail: \href{mailto:parkjhm@kunsan.ac.kr}{parkjhm@kunsan.ac.kr}}}
\affil[1]{Department of Mathematics, Kunsan National University, 558 Daehak-ro, Gunsan-si, Jeollabuk-do 54150, Republic of Korea}
\date{\small
\begin{minipage}{0.92\linewidth}\centering
This is the authors' accepted manuscript of an article published in\\
\emph{Journal of Information and Communication Convergence Engineering}, vol.~24, no.~2, pp.~155--164, June 2026.\\
Received 6 February 2026; revised 19 February 2026; accepted 27 February 2026.\\
The version of record is available at \url{https://doi.org/10.56977/jicce.2026.24.2.155}.\\[0.5ex]
This version may differ from the version of record in editorial and typographical details.\\
Please cite the published version.\\[0.5ex]
\copyright\ The Korea Institute of Information and Communication Engineering.
Distributed under the terms of the Creative Commons Attribution Non-Commercial License
(\url{http://creativecommons.org/licenses/by-nc/3.0/}), which permits unrestricted non-commercial use,
distribution, and reproduction in any medium, provided the original work is properly cited.
\end{minipage}}

\begin{document}
\maketitle

\begin{abstract}
This study proposes a Cucker--Smale model incorporating dynamic networks and the discrete $p$-Laplacian ($p \in (1,2)$) and analyzes the relationship between flocking time constraints and graph connectivity from both theoretical and experimental perspectives. We establish that this model achieves flocking within finite time and derive a quantitative relationship between the lower bound of flocking time and graph volume. Based on these results, we observe that imposing sufficiently small flocking time induces an increase in network volume, and that larger initial velocity variance requires greater volume under the same time constraint. In our experiments, we trained a transformer-based PPO policy to select network structures at each time step. Under the most stringent time constraint, the maximum average number of connections increased sharply, confirming a threshold effect. The negative correlation between convergence time and average number of connections reflects the theoretical results well. This study provides quantitative evidence for the phenomenon that time constraints induce connection mobilization.

\end{abstract}

\section*{I. Introduction}

Collective motion, commonly observed in nature such as flocking of birds
or schooling of fish, is a phenomenon in which individual agents form
velocity alignment and spatial cohesion through local interactions
without centralized control [1]. The Cucker--Smale (C-S) model,
proposed by Cucker and Smale in 2007, is a representative dynamical
model describing such flocking behavior, and its convergence
characteristics vary considerably depending on which graph mediates the
interactions among agents [2].

In practical systems, it is difficult to assume static and dense graphs
such as complete graphs due to sensing range limitations, communication
and computational costs, and environmental changes. Time-varying
interaction graphs \(\mathcal{G}(t)\), where neighborhood relationships
change over time, are more natural [3]. Consequently, the question
of ``how to select and control time-varying graphs to achieve flocking
more rapidly'' remains a central problem bridging theory and
applications.

In this context, it has been reported that in C-S type models
incorporating bonding forces, when reinforcement learning (RL) is
employed to adjust connectivity in real time, the time at which velocity
alignment error falls below a certain tolerance is advanced [10].
This prior work demonstrates that an approach treating the interaction
graph itself as a control variable can improve the ``practical
convergence time'' of flocking and serves as a direct motivation for the
present study. However, it should be noted that the flocking time in
that work is not a finite-time flocking time in the strict sense, but
rather a reaching time defined by an \(\varepsilon\)-criterion.

When the same reinforcement learning-based graph selection is applied to
the basic Cucker--Smale model without bonding forces, different behavior
emerges. The present study takes this gap as its starting point.
According to theoretical intuition, while it has not been rigorously
proven that complete connectivity induces the fastest alignment, the
expectation that more connections facilitate alignment through richer
information exchange is natural. However, the learned policies are
observed to tend to select sparse graphs maintaining only a limited
number of edges, contrary to this expectation. Indeed, in preliminary
experiments of this study, cases were observed where the
\(\varepsilon\)-criterion flocking time was delayed. This can be
interpreted as a tendency for the policy to select only the minimum
sufficient connections that satisfy the goal, given that the
optimization objective of reinforcement learning is ``achieving
\(\varepsilon\)-criterion alignment'' rather than ``approaching a
complete graph structure.'' That is, when only the ``occurrence'' of
flocking is required, the incentive to maintain large connectivity
weakens, and this choice may lead to delayed convergence time.

These observations give rise to two questions. First, why does
reinforcement learning-based graph selection not pursue complete
connectivity? Second, what learning design can guide the policy toward
increasing connectivity (edge usage)? To answer these questions, it is
necessary to first establish the theoretical relationship between time
constraints and graph connectivity.

However, the asymptotic flocking characteristics of the basic C-S model
(linear velocity coupling) create limitations. It is difficult to
naturally define the model's inherent `flocking time' as a finite value,
and the \(\varepsilon\)-criterion reaching time depends on the choice of
tolerance, making it difficult to stably quantify the time-connectivity
relationship in inequality form [2][4]. Therefore, for
theoretical analysis, this study introduces a nonlinear Cucker--Smale
model incorporating the discrete \(p\)-Laplacian on a time-varying
interaction graph \(\mathcal{G}(t)\). For \(p \in (1,2)\), finite-time
flocking holds, allowing the flocking time \(T_{f}\) to be defined at
the model level and enabling rigorous treatment of time-based
requirements. Building on this theoretical foundation, this study
proposes the hypothesis that time pressure, which directly demands
``faster flocking achievement,'' can serve as a key lever. That is, by
limiting the maximum allowed steps (\texttt{max\_steps}) and guiding learning to
satisfy alignment criteria within that limit, the policy may expand
interactions and mobilize connectivity as a resource to achieve
convergence within the limited time.

The main contributions of this study are as follows. First, in a
discrete \(p\)-Laplacian-based C-S model with \(p \in (1,2)\), we
consider time-varying interaction graphs \(\mathcal{G}(t)\) and
quantitatively derive and prove the trade-off relationship between
flocking time \(T_{f}\) and the volume of \(\mathcal{G}(t)\). This
quantitatively demonstrates that enforcing smaller flocking time
requires greater connectivity, theoretically justifying that time
pressure can lead to increased edge usage. Second, we train PPO
algorithms under various episode length conditions for nonlinear C-S
flocking, measuring average edge count and flocking achievement time.
This experimentally quantifies the effect of time pressure intensity on
graph selection and convergence characteristics. Third, we analyze how
differences in learning conditions affect convergence performance and
graph selection patterns. This presents in comparable form the
conditions under which the tendency to select sparse graphs and the
resulting convergence delay are reinforced or mitigated. Fourth, we
statistically analyze the correlation between convergence time and
average edge count. This supports at the data level that the
``convergence acceleration--edge increase'' phenomenon observed under
time pressure is not coincidental but a structural tendency.

The remainder of this paper is organized as follows. Section II reviews
the basic concepts of nonlinear C-S models, graph theory, and core
elements of reinforcement learning. Section III defines the dynamic
network-based nonlinear Cucker-Smale model and presents the theoretical
relationship between flocking time and graph volume. Section IV presents
comparative experimental results according to learning \texttt{max\_steps}
settings, and Section V summarizes the research and discusses future
directions.

\section*{II. Background}

This section reviews the dynamic models, graph notation, and minimal
concepts of reinforcement learning used in subsequent sections.

\subsection*{A. The C-S Model and Nonlinear Extensions}

Consider \(N\) agents moving in \(d\)-dimensional space. Let the
position and velocity of the \(i\)-th agent (\(i = 1,\ldots,N\)) be
\(x_{i} = (x_{i}^{1},\ldots,x_{i}^{d})^{T} \in \mathbb{R}^{d}\) and
\(v_{i} = (v_{i}^{1},\ldots,v_{i}^{d})^{T} \in \mathbb{R}^{d}\),
respectively. For time \(t \geq 0\), the C-S model is given as follows

\[x_{i}' = v_{i}\]

\[v_{i}' = \sum_{j = 1}^{N}\psi( \|x_{j} - x_{i}\| )(v_{j} - v_{i})\]

Here, \(x_{i}'\) and \(v_{i}'\) denote the time derivatives of the
position and velocity vectors, respectively,
\(\|\cdot\|\) is the standard Euclidean distance, and
the function \(\psi\) is the interaction function defined as

\[\psi(r) = \frac{1}{(1 + r^{2})^{\frac{\beta}{2}}},\quad\beta > 0.\]

The parameter \(\beta\) controls the decay rate of the interaction;
larger \(\beta\) results in more rapid decrease of influence with
distance. Under such interaction, the system is said to exhibit flocking
when the positions and velocities of the agents satisfy the following
two conditions

\[\lim_{t \to \infty}\max_{1 \leq i,j \leq N} \|v_{i} - v_{j}\| = 0\]

\[\sup_{t \geq 0}\max_{1 \leq i,j \leq N} \|x_{i} - x_{j}\| \leq M < \infty\]

Cucker and Smale showed that when the interaction decay rate is
\(\beta \in (0,1)\), flocking occurs regardless of initial conditions,
and when \(\beta \geq 1\), flocking is proven to occur under certain
initial conditions [2]. Other topics related to the C-S model
include mathematical analysis of flocking [11], multi-cluster
flocking [12], formation control [13], collision avoidance
[14], among various other studies [15], [16].

In the basic C-S model, the rate of velocity change of each agent is
exactly proportional to the velocity difference from its neighbors.
However, in real collective motion, the response intensity may vary
nonlinearly depending on the magnitude of the velocity difference. For
example, saturation effects may occur when the velocity difference is
very large, while stronger responses may be needed for fine adjustments
when it is very small. Kim and Park [4] proposed a model introducing
a nonlinear weight \(|v_{j}^{k} - v_{i}^{k}|^{p - 2}\) on the velocity
difference to reflect such nonlinear effects. This model is closely
related to the structure of the discrete \(p\)-Laplacian operator, and
flocking occurrence and convergence speed are determined by the value of
\(p\). The nonlinear C-S model is defined as follows

\[\left( x_{i}^{k} \right)^{'} = v_{i}^{k},\ \ \ (v_{i}^{k})' = \sum_{j = 1}^{N}\psi_{ij}|v_{j}^{k} - v_{i}^{k}|^{p - 2}(v_{j}^{k} - v_{i}^{k}),\]

where \(\psi_{ij} = \psi( \|x_{j} - x_{i}\| )\).

\subsection*{B. Graph Laplacian and Connectivity Measures}

This section introduces the graph Laplacian, which quantifies the
interaction graph.

The adjacency matrix \(A\) is defined as an \(N \times N\) matrix with
entry \(1\) (or a weight value) if an edge exists between vertices \(i\)
and \(j\), and \(0\) otherwise. The degree matrix \(D\) is a diagonal
matrix whose diagonal entries are the degrees of each vertex, i.e.,
\(d_{i} = \sum_{j = 1}^{N}A_{ij}\). The volume of graph \(\mathcal{G}\)
is defined as the sum of degrees of all vertices,
\(\operatorname{vol}(\mathcal{G}) = \sum_{i = 1}^{N}d_{i}\). The graph Laplacian
matrix \(L\) is \(L = D - A\).

If a real number \(\lambda\) and a nonzero vector
\(u = (u_{1},\ldots,u_{N})^{T} \in \mathbb{R}^{N}\) satisfy

\begin{equation*}\sum_{j = 1}^{N}A_{ij}\left( u_{i} - u_{j} \right) = \lambda u_{i}\tag{2.1}\end{equation*}

for \(i = 1,\ldots,N\), then \(\lambda\) is called an eigenvalue of
\(L\) and \(u\) is called an eigenvector. This can be written in matrix
form as \(Lu = \lambda u\). Since the graph Laplacian \(L\) is a
symmetric matrix, all eigenvalues are real and nonnegative. Arranged in
increasing order

\[0 = \lambda_{1} \leq \lambda_{2} \leq \cdots \leq \lambda_{N - 1} \leq \lambda_{N}\]

For a complete graph, \(\lambda_{2} = \cdots = \lambda_{N} = N\), and
for a connected graph, \(\lambda_{2} > 0\). In general,
\(\lambda_{N} \leq N\) holds [5].

By equation (2.1), the eigenvalues \(\lambda_{2}\) and \(\lambda_{N}\)
are expressed as the following Rayleigh quotients

\[
\lambda_{2} \equiv \min_{u \in \mathbf{1}^{\perp}}\frac{u^{T}Lu}{u^{T}u} = \min_{u \in \mathbf{1}^{\perp}}\frac{\sum_{i < j}A_{ij}(u_{i} - u_{j})^{2}}{\sum_{i}u_{i}^{2}},
\]

\[
\lambda_{N} \equiv \max_{u \in \mathbf{1}^{\perp}}\frac{u^{T}Lu}{u^{T}u} = \max_{u \in \mathbf{1}^{\perp}}\frac{\sum_{i < j}A_{ij}(u_{i} - u_{j})^{2}}{\sum_{i}u_{i}^{2}},
\]

where \(\mathbf{1} = (1,\ldots,1)^{T} \in \mathbb{R}^{N}\). Therefore, for all
\(u \in \mathbf{1}^{\perp}\),

\[
\lambda_{2} \|u\|^{2} \leq u^{T}Lu \leq \lambda_{N} \|u\|^{2}.
\]

For further details on graph theory, see [5].

\subsection*{C. Reinforcement Learning Algorithms and Policy Networks}

The problem of selecting time-varying graphs can be viewed as a dynamic
decision-making problem, which this study addresses using reinforcement
learning. Reinforcement learning is a process of learning a policy
\(\pi\) that maximizes cumulative reward through interaction with the
environment, enabling adaptive responses to dynamic conditions such as
communication constraints and environmental changes.

The reinforcement learning problem is formalized as a Markov Decision
Process (MDP) \((\mathcal{S},\mathcal{A},P,R,\gamma)\). Here, \(\mathcal{S}\) is
the state space, \(\mathcal{A}\) is the action space, \(P(s'|s,a)\) is
the transition probability, \(R\) is the immediate reward, and
\(\gamma \in [ 0,1)\) is the discount factor. At each time \(t\),
the agent observes state \(s_{t}\in\mathcal{S}\), selects action
\(a_{t}\in\mathcal{A}\), receives reward \(r_{t} = R(s_{t},a_{t})\),
and transitions to a new state \(s_{t + 1}\).

To quantify policy performance, we use the state-value function
\(V^{\pi}\) and the action-value function \(Q^{\pi}\). \(V^{\pi}\) is
the expected cumulative reward when following policy \(\pi\) from state
\(s\), and \(Q^{\pi}\) is the expected cumulative reward when taking a
specific action \(a\) in state \(s\) and then following \(\pi\) [6].

PPO, used in this study, is an actor-critic-based policy gradient method
that stabilizes learning by suppressing excessive policy updates. The
actor approximates the policy \(\pi_{\theta}\), and the critic
approximates the value function \(V_{\phi}\). The critic is gradually
improved through Temporal Difference (TD) learning, and the 1-step TD
error is

\[\delta_{t} = r_{t} + \gamma V_{\phi}(s_{t + 1}) - V_{\phi}(s_{t})\]

The TD error \(\delta_{t}\) represents the difference between the
current value estimate \(V_{\phi}(s_{t})\) approximated by value
function parameters \(\phi\) and \(r_{t} + \gamma V_{\phi}(s_{t + 1})\)
based on actual rewards [6].

Generalized Advantage Estimation (GAE) is used for estimating the
advantage function \(\widehat{A}_{t}\) for policy improvement. This
collects TD errors over multiple steps and averages them with greater
weight on closer time steps

\[\widehat{A}_{t} = \delta_{t} + (\gamma\lambda)\widehat{A}_{t + 1}\]

This recursive formula adds to the TD error \(\delta_{t}\) at current
time \(t\) and future advantages reduced by the discount rate
(\(\gamma\)) and weight (\(\lambda\)).

The return value is defined as

\[\widehat{R}_{t} = \widehat{A}_{t} + V_{\phi}(s_{t})\]

This serves as the target value for critic learning [6][7].

Instead of the standard policy objective function, PPO optimizes the
following ratio-based loss function

\[L^{CLIP}(\theta) = \mathbb{E}_{t}\left[ \min\left( \rho_{t}(\theta)\widehat{A}_{t},\operatorname{clip}(\rho_{t}(\theta),1 - \epsilon,1 + \epsilon)\widehat{A}_{t} \right) \right]\]

where
\(\rho_{t}(\theta) = \frac{\pi_{\theta}(a_{t}|s_{t})}{\pi_{\theta_{\mathrm{old}}}(a_{t}|s_{t})}\)
is the probability ratio between the new and old policies. \(\epsilon\)
is a small positive number that limits the clipping range of policy
change. That is, when \(\rho_{t}(\theta)\) exceeds the
\(1 \pm \epsilon\) range, the magnitude of policy update is suppressed
(clipping), stabilizing learning [7].

For approximating the policy \(\pi_{\theta}\), this study uses a
transformer-based encoder as the network architecture. The transformer
can flexibly reflect dependencies among agents through the attention
mechanism, making it suitable for the present problem where a fixed
adjacency structure is not assumed (time-varying \(\mathcal{G}(t)\))
[8][9].

\section*{III. Dynamic Network-Based Nonlinear Cucker--Smale Model}

This section introduces an extended model based on the nonlinear C-S
model proposed by Park and Kim [4], considering time-varying network
connection structures. While the original study [4] assumed a
complete graph, the present study introduces a time-dependent adjacency
matrix \(A_{ij}(t)\) to handle situations where connection structures
change over time through reinforcement learning.

Therefore, in this study, we assume that for \(N\) agents, the positions
\(x_{i}\) and velocities \(v_{i}\) satisfy

\begin{equation*}(x_{i}^{k})' = v_{i}^{k},\tag{3.1}\end{equation*}

\begin{equation*}(v_{i}^{k})' = \sum_{j = 1}^{N}{A_{ij}\psi_{ij}}|v_{j}^{k} - v_{i}^{k}|^{p - 2}\left( v_{j}^{k} - v_{i}^{k} \right).\tag{3.2}\end{equation*}

Here, we assume that the graph obtained from the adjacency matrix is
connected for all time \(t\). This corresponds to the condition of a
connected graph defined in Section II, and under this assumption, the
second eigenvalue \(\lambda_{2}(t) > 0\) of the graph Laplacian \(L\) is
guaranteed for all time \(t \geq 0\) [5]. Furthermore, we restrict
the parameter ranges to \(p \in (1,2)\) and \(\beta \in (0,1]\),
within which finite-time flocking occurs regardless of initial
conditions (Theorem 3.1). Moreover, a lower bound on the flocking time
has a quantitative relationship with the graph volume (Theorem 3.2),
providing theoretical justification for the phenomenon where time
pressure induces increased connectivity in reinforcement learning
experiments.

For convenience of analysis, we introduce fluctuation variables. Let the
mean position and mean velocity be
\(\bar{x} = \frac{1}{N}\sum_{j = 1}^{N}x_{j}\) and
\(\bar{v} = \frac{1}{N}\sum_{j = 1}^{N}v_{j}\), respectively, and define
the fluctuations from the mean as

\[\widehat{x}_{i} = x_{i} - \bar{x},\quad\widehat{v}_{i} = v_{i} - \bar{v}.\]

In particular, using the symmetry of the adjacency matrix in equation
(3.2) and the odd function property of the interaction term, one can
verify that the mean velocity \(\bar{v}\) is conserved in time. From
this, the fluctuation variables satisfy the following system

\begin{equation*}(\widehat{x}_{ik})' = \ \widehat{v}_{ik}\tag{3.3}\end{equation*}

\begin{equation*}(\widehat{v}_{ik})' = \sum_{j = 1}^{N}{A_{ij}\psi_{ij}\left| \widehat{v}_{jk} - \widehat{v}_{ik} \right|^{p - 2}\left( \widehat{v}_{jk} - \widehat{v}_{ik} \right)}.\tag{3.4}\end{equation*}

We now state the main results concerning flocking.

\medskip\noindent\textbf{Theorem 3.1} (\emph{Finite-Time Flocking}). Suppose \(\beta \in (0,1]\),
\(p \in (1,2)\), and the graph is connected for all \(t \geq 0\). Then
velocity alignment is achieved within finite time. More precisely, for
\(t \geq T^{*}\), we have
\(\|\widehat{v}(t)\|_{2} = 0\), where the upper bound
\(T^{*}\) on the flocking time is

\[T^{*} = \frac{2N^{p} \|\widehat{v}(0)\|_{2}^{2 - p}}{(2 - p)C_{*}}\]

Here, the velocity fluctuation energy is defined as
\(\|\widehat{v}(t)\|_{2}^{2}:= \sum_{i = 1}^{N} \|\widehat{v}_{i}(t)\|_{2}^{2}\),
and \(C_{*}\) is given by

\[
C_{*} = \begin{cases}
\left( \psi_{\min}(0)^{\frac{\beta - 1}{\beta}} + \frac{2\sqrt{2}(1 - \beta)N^{p}}{3 - p} \|\widehat{v}(0)\|_{2}^{3 - p} \right)^{\frac{- \beta}{1 - \beta}} & \text{if }0 < \beta < 1, \\
\psi_{\min}(0)\exp\left( - \frac{2\sqrt{2}N^{p}}{3 - p} \|\widehat{v}(0)\|_{2}^{3 - p} \right) & \text{if }\beta = 1.
\end{cases}
\]

\noindent\emph{Proof.} The work [4] addressed the case of a complete graph
(i.e., \(A_{ij}(t) = 1\) for all \(i,j = 1,\ldots,N\) and \(t \geq 0\)),
in which case \(\lambda_{2} = N > 0\). In the present paper,
\(\lambda_{2}\) depends on time, but since we assume the graph is
connected for all time, applying the same method as in the proof of
[5, Lemma 1.9], one can show that
\(\lambda_{2}(t) \geq \frac{1}{2N^{2}}\) for all \(t \geq 0\). Since the
positivity of \(\lambda_{2}\) is secured, applying the remaining logic
of [4] identically yields this theorem. $\blacksquare$

\medskip\noindent\textbf{Theorem 3.2} (\emph{Lower Bound on Flocking Time}). Suppose
\(\beta \in (0,1]\), \(p \in (1,2)\), and the graph is connected
for all \(t \geq 0\). Then there exists a lower bound \(T_{*} > 0\) on
the flocking time satisfying

\[\int_{0}^{T_{*}}\operatorname{vol}(\mathcal{G}(s))^{1 - p/2}ds = \frac{2 \|\widehat{v}(0)\|_{2}^{2 - p}}{(2 - p)d^{1 - p/2}(2N)^{p/2}}.\]

\noindent\emph{Proof.} Taking the inner product of both sides of equation (3.4)
with \(\widehat{v}_{i}\) and summing over \(i\), we obtain the
following energy equation

\[
\frac{d}{dt} \|\widehat{v}\|_{2}^{2} = - \sum_{i,j}A_{ij}\psi_{ij} \|\widehat{v}_{j} - \widehat{v}_{i}\|_{p}^{p}.
\]

Here, \(\|u\|_{p}^{p}:= \sum_{k = 1}^{d}|u_{k}|^{p}\).
Since \(\psi_{ij} \leq 1\),

\[
\frac{d}{dt} \|\widehat{v}\|_{2}^{2} \geq - \sum_{i,j}A_{ij} \|\widehat{v}_{j} - \widehat{v}_{i}\|_{p}^{p}.
\]

For \(1 \leq p \leq 2\), we have
\(\|a\|_{p} \leq d^{1/p - 1/2} \|a\|_{2}\)
for any \(a \in \mathbb{R}^{d}\), so

\[
\frac{d}{dt} \|\widehat{v}\|_{2}^{2} \geq - d^{1 - p/2}\sum_{i,j}A_{ij} \|\widehat{v}_{j} - \widehat{v}_{i}\|_{2}^{p}.
\]

By Jensen's inequality,

\[
\sum_{i,j}A_{ij} \|\widehat{v}_{j} - \widehat{v}_{i}\|_{2}^{p} \leq \operatorname{vol}(\mathcal{G}(t))^{1 - \frac{p}{2}}\left( \sum_{i,j}A_{ij} \|\widehat{v}_{j} - \widehat{v}_{i}\|_{2}^{2} \right)^{\frac{p}{2}}.
\]

By the definition of fluctuation, \(\sum_{i}\widehat{v}_{ik} = 0\),
so for each \(k\), the vector
\((\widehat{v}_{1k},\ldots,\widehat{v}_{Nk})^{T} \in \mathbf{1}^{\perp}\).
Therefore, by the Rayleigh quotient in Section II.B,
\(\lambda_{N} \leq N\), hence

\[
\sum_{i,j}A_{ij} \|\widehat{v}_{j} - \widehat{v}_{i}\|_{2}^{2} \leq 2N \|\widehat{v}\|_{2}^{2}.
\]

Therefore

\[
\frac{d}{dt} \|\widehat{v}\|_{2}^{2} \geq - d^{1 - p/2}(2N)^{p/2}\operatorname{vol}(\mathcal{G}(t))^{1 - p/2} \|\widehat{v}\|_{2}^{p}.
\]

Dividing both sides by \(\|\widehat{v}\|_{2}^{p}\) and
integrating from 0 to \(t\),

\[
\|\widehat{v}(t)\|_{2}^{2 - p} \geq \|\widehat{v}(0)\|_{2}^{2 - p} - C\int_{0}^{t}\operatorname{vol}(\mathcal{G}(s))^{1 - \frac{p}{2}}ds,
\]

where \(C = \frac{(2 - p)d^{1 - \frac{p}{2}}(2N)^{\frac{p}{2}}}{2}\).
Since the graph is connected,
\(\operatorname{vol}(\mathcal{G}(t)) \geq 2(N - 1) > 0\), and as
\(t \rightarrow \infty\),
\(\int_{0}^{t}\operatorname{vol}(\mathcal{G}(s))^{1 - p/2}ds \rightarrow \infty\).
Therefore, there exists a \(T_{*}\) at which the right-hand side becomes
0, and this serves as a lower bound on the finite flocking time \(T\)
that exists by Theorem 3.1. $\blacksquare$

\medskip\noindent\textbf{Remark 3.3} (\emph{Implications of the Time-Volume Trade-off}). Theorem 3.2
reveals a trade-off relationship between flocking time and graph volume.
As mentioned in the introduction, in reinforcement learning we impose
time pressure by limiting the maximum allowed steps. Therefore, letting
\(T_{\max}\) denote the allowed time until flocking, Theorem 3.2 has the
following implications.

(i) The role of \(T_{*}\) as a threshold: Let the actual flocking time
be \(T\), so \(T \geq T_{*}\). Given a time limit \(T_{\max}\), if
\(T_{\max} > T_{*}\), the lower bound condition is satisfied, so volume
increase is not necessarily required. On the other hand, if
\(T_{\max} < T_{*}\), flocking cannot occur within the given time with
the current graph structure. In this case, to satisfy the equation while
reducing \(T_{*}\) to below \(T_{\max}\), the integral of
\(\operatorname{vol}(G(s))^{1 - p/2}\) must remain the same, so the volume must
increase. Therefore, moderately reducing the time limit has little
effect on volume increase; it must be reduced sufficiently below
\(T_{*}\) to force volume increase.

(ii) Effect of initial velocity variance: The right-hand side contains
\(\|\widehat{v}(0)\|_{2}^{2 - p}\), and since
\(p \in (1,2)\), we have \(2 - p > 0\). Therefore, larger initial
velocity variance \(\|\widehat{v}(0)\|_{2}\) increases
the right-hand side of the equation, requiring greater volume under the
same time limit.

(iii) Connection to reinforcement learning experiments: Combining (i)
and (ii), when training graph structures through reinforcement learning
in an environment with large initial velocity variance and strict time
limits, an increase in volume (number of edges) is predicted. This
prediction will be verified through experiments in Section IV.

\section*{IV. Experiments and Results}

\subsection*{A. Experimental Setup}

This section presents the specific settings and results of experiments
conducted by implementing the dynamic network-based nonlinear C-S model
from Section III as a reinforcement learning environment. The main
objective of the experiments is to quantitatively analyze the effect of
episode length constraints (\texttt{max\_steps}) imposed during training on the
graph selection patterns and flocking convergence time of the learned
policies.

In particular, shorter \texttt{max\_steps} act as `time pressure' requiring the
same convergence goal to be achieved more quickly, and we verify whether
policies mobilize connectivity as a resource for convergence. We use
settings where finite-time flocking holds in the nonlinear C-S model to
rigorously define flocking time and use the theoretical observation that
the lower bound on flocking time has an inverse relationship with graph
connectivity as the basis for interpreting results.

The simulation assumes a dynamic network where the adjacency matrix
changes over time. The number of agents is set to \(N = 20\), and
initial positions are randomly generated from a uniform distribution in
the range \([ - 10,10] \times [ - 10,10]\).
Initial velocities \(v_{i}(0)\) are generated using speeds randomly
selected from \([ 0,10]\) and angles uniformly distributed
in \([ 0,2\pi]\). The time step is fixed at
\(\Delta t = 0.1\).

The observation space is defined as an \(N \times 4\) matrix containing
position and velocity information of each agent. For each agent \(i\), a
4-dimensional vector \((x_{i},y_{i},v_{x,i},v_{y,i})\) is provided as
observation. The action space is defined as \texttt{MultiBinary(n\_edges)}, where
\(\texttt{n\_edges} = \frac{N(N - 1)}{2} = 190\). Each action is
represented as a 190-dimensional binary vector, where each element
indicates the presence (1: connected, 0: not connected) of an undirected
edge. This action vector is internally converted to a symmetric
adjacency matrix \(A(t)\) and used for the physics update in equation
(3.2).

Each episode proceeds for up to \texttt{max\_steps}, and in this study,
experiments are conducted with various settings
\(\texttt{max\_steps} \in \left\{ 10,20,30,50,100,1000 \right\}.\) An
episode terminates when one of two conditions is met. First, if velocity
variance decreases below a threshold achieving flocking, the episode
terminates as success. Second, if flocking is not achieved by
\texttt{max\_steps}, the episode terminates as failure.

\subsection*{B. Reward Function Design}

The reward function is designed so that agents clearly recognize the
learning objective and develop effective policies. The total reward is
defined differently for success and failure as follows

\noindent On success:
\[
r_{\mathrm{total}}(t) = r_{\mathrm{align}}(t) + r_{\mathrm{success}}(t),
\]

\noindent On failure:
\[
r_{\mathrm{total}}(t) = r_{\mathrm{failure}}(t).
\]

First, the alignment reward \(r_{\mathrm{align}}\) is defined as

\[
r_{\mathrm{align}}(t) = w_{\mathrm{align}} \times 10.0 \times [ V(t - 1) - V(t)].
\]

Here, \(w_{\mathrm{align}} = 5.0\) is the weight for alignment reward, and 10.0
is a scaling factor. \(V = \operatorname{Var}(v_{x}) + \operatorname{Var}(v_{y})\) is the sum of
velocity variances along each axis, where the variance for each axis is
defined as
\(\operatorname{Var}(v_{x}) = \frac{1}{N}\sum_{i = 1}^{N}(v_{x,i} - \bar{v}_{x})^{2}\).
This reward structure provides positive reward when acting to decrease
velocity variance and negative reward when acting to increase it,
guiding agents to continuously learn toward flocking achievement.

At episode termination, success reward or failure penalty is assigned.
If flocking is achieved with velocity variance decreasing below
threshold \(\texttt{vel\_threshold} = 0.1\), a bonus of
\(r_{\mathrm{success}} = 100.0\) is given. If flocking is not achieved by
\texttt{max\_steps}, a penalty of \(r_{\mathrm{failure}} = - 100.0\) is assigned without
alignment reward. This reward structure enables agents to learn the time
pressure of achieving flocking within a limited time.

\subsection*{C. Hyperparameter Settings}

This study conducts experiments with nonlinear parameter \(p = 1.5\) and
decay coefficient \(\beta = 0.5\). This satisfies the conditions of
Theorem 3.1, ensuring that experiments are conducted in the region where
finite-time flocking is theoretically guaranteed. In particular, in this
range, boundedness of distance is also secured within finite time at the
moment velocity alignment is achieved, which is the theoretical basis
for using only the velocity alignment metric in the reward function of
Section IV.B.

Such settings are also important from an experimental design
perspective. By theoretically excluding potential confounding factors
such as distance boundedness issues and sensitivity to initial
conditions, observed performance differences are primarily explained by
differences in episode length constraints (\texttt{max\_steps}) imposed during
learning, rather than by chance in initialization or instability of
dynamics. That is, we can more clearly compare how graph selections of
policies learned at different \texttt{max\_steps} respond to ``given time
constraints.''

In our experiments, connectivity is not enforced as an explicit action
constraint. Theorem 3.2 derives the relationship between time
constraints and graph volume under the connectivity assumption,
providing theoretical grounds for the hypothesis that time pressure
induces increased connectivity. To test this hypothesis, the policy
freely selects adjacency matrices that may include disconnected graphs,
and episodes in which flocking was achieved are used for analysis
regardless of connectivity. Since this setting is applied identically
across all \texttt{max\_steps} conditions, differences between conditions arise
from differences in the intensity of time pressure.

Learning and implementation settings are fixed as follows. PPO
algorithms were trained independently for each of the 6 \texttt{max\_steps}
conditions defined in Section IV.A. For \texttt{max\_steps} of 10, 20, 30, and
100, learning was conducted for 270, 400, 380 and 300 iterations,
respectively, while 200 iterations were used for the remaining
conditions. For conditions with smaller \texttt{max\_steps}, the initial failure
rate is high, resulting in slower convergence, so the number of
iterations was increased to ensure sufficient learning. During training,
models with the highest episode reward mean were selected, and these
models were used in the evaluation phase.

In the test phase, to ensure fair comparison, all learned policies were
evaluated in the same test environment. The test environment \texttt{max\_steps}
was fixed at 100. For each policy, episodes were run from 5000 different
initial conditions (seeds) to collect statistically meaningful results.

\begin{table}[ht]
\centering
\caption*{Table I. Parameter setting}
\begin{tabular}{ll}
\toprule
Parameter & Value \\
\midrule
PPO learning rate & 0.0003 \\
Batch size & 4000 \\
Discount factor & 0.99 \\
Gradient clipping & 0.5 \\
Training iterations & 200 or condition-specific \\
\bottomrule
\end{tabular}
\end{table}

\subsection*{D. Experimental Results and Analysis}

The test performance of policies learned under each \texttt{max\_steps} condition
is summarized in Table II. Success rate represents the percentage of
episodes achieving flocking out of 5000 test episodes.

\begin{table}[ht]
\centering
\caption*{Table II. Test results by training \texttt{max\_steps}. \texttt{avg\_edges} is the time-averaged number of connections selected at each step within an episode, and \texttt{avg\_edges\_min} and \texttt{avg\_edges\_max} are the minimum and maximum values of \texttt{avg\_edges} observed across 5000 test episodes.}
\begin{tabular}{cccc}
\toprule
Training \texttt{max\_steps} & Success rate (\%) & \texttt{avg\_edges\_min} & \texttt{avg\_edges\_max} \\
\midrule
10 & 98.76 & 50.39 & 177.85 \\
20 & 90.58 & 34.01 & 149.23 \\
30 & 96.08 & 33.25 & 146.12 \\
50 & 99.86 & 81.29 & 154.66 \\
100 & 97.64 & 33.51 & 119.67 \\
1000 & 97.96 & 30.31 & 149.18 \\
\bottomrule
\end{tabular}
\end{table}

A notable point in Table II is the change in \texttt{avg\_edges\_max}. Even when
reducing \texttt{max\_steps} from 1000 to 20, \texttt{avg\_edges\_max} shows little
variation in the range of 141 to 155, but at \texttt{max\_steps}=10, it sharply
increases to 177.85. This demonstrates, as predicted in Remark 3.3 (i),
that moderately reducing the time limit has little effect on volume
increase, and volume increase is induced only when reduced sufficiently
below the threshold \(T_{*}\).

However, \texttt{avg\_edges\_min} does not show a consistent increase pattern
with decreasing \texttt{max\_steps}. This means that volume does not increase in
all episodes, but only in some. As mentioned in Remark 3.3 (ii), larger
initial velocity variance requires greater volume, so volume increase is
prominent in episodes with difficult initial conditions, while sparse
connections may suffice for convergence in episodes with easy initial
conditions.

Next, we analyze the relationship between number of connections and
convergence time. Fig. 1 presents scatter plots of (\texttt{terminated\_step},
\texttt{avg\_edges}) pairs observed under each learning condition. Here,
\texttt{terminated\_step} is the convergence step at which the episode
terminated.

The dashed lines at the top of Fig. 1 indicate baselines. full edges=190
is the number of connections in a complete graph, and full step
min/mean/max=4.0/6.5/9.0 is the range of convergence steps when using a
complete graph at every step. That is, using maximum connections
achieves convergence within 4 to 9 steps.

\begin{figure}[ht]
\centering
\includegraphics[width=0.75\linewidth]{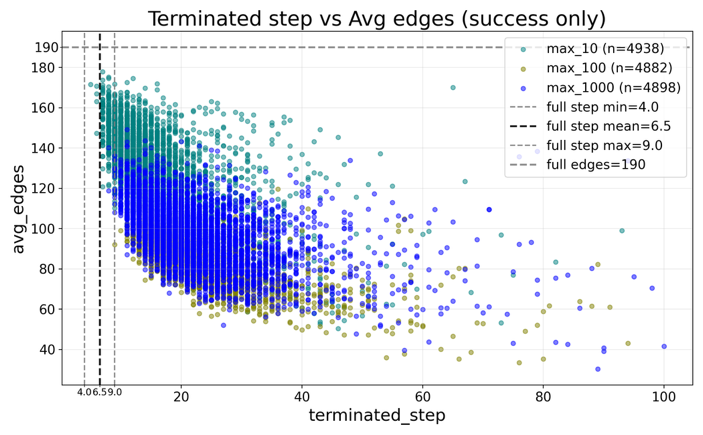}
\caption*{Fig. 1. Scatter plot of test results for \texttt{max\_steps} 10, 100, 1000}
\end{figure}

A negative correlation is observed across all three conditions in the
scatter plot. Episodes with fast convergence (small \texttt{terminated\_step})
have higher \texttt{avg\_edges}, and \texttt{avg\_edges} generally decreases as
\texttt{terminated\_step} increases. The distribution regions also differ by
condition. \texttt{max\_10} (cyan) is distributed in the upper region with
\texttt{terminated\_step} 4 to 30 and \texttt{avg\_edges} 100 to 190, while \texttt{max\_1000}
(blue) is distributed in the middle region with \texttt{terminated\_step} 10 to
50 and \texttt{avg\_edges} 50 to 150.

\begin{figure}[ht]
\centering
\includegraphics[width=0.75\linewidth]{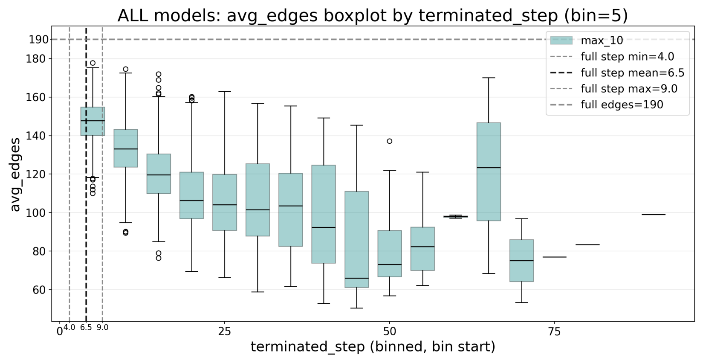}
\caption*{Fig. 2. Box plot for \texttt{max\_steps} 10}
\end{figure}

To quantitatively confirm this trend, Pearson correlation coefficients
between \texttt{terminated\_step} and \texttt{avg\_edges} were calculated for successful
episodes across the 6 learning conditions. A negative correlation was
observed in all conditions. In more detail, \texttt{max\_steps}=10 and 100 showed
relatively stronger correlations of $\rho=-0.59$ and $\rho=-0.68$, respectively,
while \texttt{max\_steps}=1000 showed $\rho=-0.49$. The remaining conditions,
\texttt{max\_steps}=20, 30, and 50, showed relatively weaker correlations of
$\rho=-0.39$, $\rho=-0.30$, and $\rho=-0.37$, respectively. Despite the variation in
magnitude, this consistent negative correlation is consistent with the
prediction of Theorem 3.2.

While scatter plots intuitively show overall correlation trends, it is
difficult to examine distributions by interval. To supplement this,
Figs. 2--3 bin \texttt{terminated\_step} with bin size 5 and present box plots of
\texttt{avg\_edges} distributions for episodes in each bin.

\begin{figure}[ht]
\centering
\includegraphics[width=0.75\linewidth]{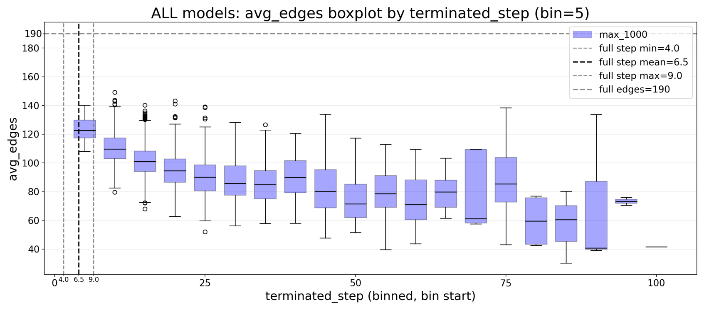}
\caption*{Fig. 3. Box plot for \texttt{max\_steps} 1000}
\end{figure}

A common pattern of increasing \texttt{avg\_edges} median with decreasing
\texttt{terminated\_step} is observed in both figures. However, distribution
regions differ by condition. Fig. 2 (\texttt{max\_steps}=10) is distributed
overall in a high \texttt{avg\_edges} region, while Fig. 3 (\texttt{max\_steps}=1000) is
distributed in a lower region. There are also differences in connection
usage even within the same convergence interval. For example, in the
early interval (\texttt{terminated\_step} 0 to 10), the median of \texttt{max\_10} is
approximately 145, while that of \texttt{max\_1000} is approximately 120. This
shows that policies learned under \texttt{max\_10} use more connections even when
achieving the same fast convergence.

Synthesizing these analysis results, the connection usage tendency of
policies varies according to time constraints imposed during training.
Policies learned under strong time constraints (\texttt{max\_steps}=10) actively
use connection resources to shorten convergence time, while policies
learned under weak time constraints (\texttt{max\_steps}=1000) converge in a
direction that conserves connections within the range of maintaining
success. Such connection usage tendencies are observable even in the
same test environment.

\section*{V. Conclusion and Discussion}

This study analyzed the relationship between time constraints and graph
connectivity in the nonlinear Cucker--Smale model from both theoretical
and experimental perspectives.

Theoretically, Theorem 3.2 presented an integral equation between the
lower bound on flocking time and graph volume. Based on this, Remark 3.3
derived two predictions: (i) moderately reducing the time limit has
little effect on volume increase, and volume increase is induced only
when reduced sufficiently below the threshold \(T_{*}\); (ii) larger
initial velocity variance requires greater volume under the same time
limit.

The experimental results support these theoretical predictions. In Table
II, even when reducing \texttt{max\_steps} from 1000 to 20, \texttt{avg\_edges\_max}
showed little variation, but at \texttt{max\_steps}=10 it sharply increased,
confirming the threshold effect of Remark 3.3 (i). The lack of a
consistent increase pattern in \texttt{avg\_edges\_min} is consistent with Remark
3.3 (ii). The negative correlation between \texttt{terminated\_step} and
\texttt{avg\_edges} observed in the scatter plot of Fig. 1 is consistent with the
results of Theorem 3.2.

These results support the explanatory framework that ``time constraints
induce connection mobilization.'' Policies learned under strong time
constraints actively use connection resources to shorten convergence
time, while policies learned under weak time constraints converge in a
direction that conserves connections. This study is significant in that
it presents such phenomena quantitatively together with theoretical
justification. Furthermore, while existing Cucker-Smale studies have
primarily focused on ``does flocking occur?'', this study presents a
resource allocation perspective asking, ``how many connections are
needed for flocking?'' If connections are viewed as cost, the
time-connection trade-off may have implications for communication cost
optimization in multi-agent systems such as drone swarms or robot
swarms. Moreover, reinforcement learning policies did not explicitly
learn Theorem 3.2, and the connectivity assumption was not enforced in
the experiments. Nevertheless, trends consistent with the theoretical
predictions were observed, suggesting that the relationship between time
constraints and connectivity increase may hold even under relaxed
connectivity assumptions.

Future work may extend to analyzing structural characteristics such as
\(\lambda_{2}\) and degree distribution beyond \texttt{avg\_edges} or explicitly
introducing connection costs into the reward to directly optimize the
trade-off between time and connection resources.

\section*{Acknowledgements}

This work was supported by National Research Foundation of Korea (NRF)
grant funded by the Korean government (MSIT) (No. 2021R1A2C1093929).

\end{document}